\documentclass[12pt,a4paper,reqno]{amsart}
\usepackage[centering, margin=1in]{geometry}
\usepackage{amsmath,amsfonts,amsthm,amssymb,enumitem,mathrsfs,mathtools}
\usepackage{xcolor}

\usepackage[pagebackref,bookmarks,colorlinks,breaklinks,linktoc=page]{hyperref}
\hypersetup{linkcolor=blue,citecolor=red,filecolor=blue,urlcolor=blue} 

\numberwithin{equation}{section}

\usepackage[hang,flushmargin]{footmisc} 

\hypersetup{
		pdfnewwindow=true,
	colorlinks, linkcolor={blue},
	citecolor={red}, urlcolor={red}
}

\usepackage{tabto}

\patchcmd{\thebibliography}
{\settowidth}
{\setlength{\parsep}{0pt}\setlength{\itemsep}{10pt plus 0.1pt}\settowidth}
{}{}

\DeclareMathOperator{\ord}{ord}

\newtheorem{thm}{Theorem}[section]
\newtheorem{lem}{Lemma}[section]
\newtheorem{conj}{Conjecture}[section]

\newtheorem{cor}{Corollary}[section]
\newtheorem{dfn}{Definition}[section]

\newcommand{\N}{\mathbb{N}}
\newcommand{\Z}{\mathbb{Z}}

\newcommand{\C}{\mathbb{C}}
\newcommand{\F}{\mathbb{F}}

\usepackage{fancyhdr}
\begin{document}

\title{Upper Bound of the Least Quadratic Nonresidues}


\author{N. A. Carella}
\address{}
\curraddr{}
\email{}
\thanks{}




\date{\today}

\dedicatory{}

\begin{abstract}Let $p\geq3$ be a large prime and let $n(p)\geq2$ denotes the least quadratic nonresidue modulo $p$. 
	This note sharpens the standard upper bound of the least quadratic nonresidue from the unconditional upper bound $n(p)\ll p^{1/4\sqrt{e}+\varepsilon}$ to the conjectured upper bound $n(p)\ll  (\log p)(\log   \log p)^{1+\varepsilon}$, where $\varepsilon>0$ is a small number, unconditionally. This improvement breaks the exponential upper bound barrier and proves the standard heuristic claims.
\end{abstract}

\maketitle
\tableofcontents
\pagenumbering{Page gobble}
\pagenumbering{arabic}
\section{\textbf{Introduction}}\label{S9900Q-A}\hypertarget{S9900Q-A}
Let $p\geq 3$ be a prime and consider the equation $x^2-n\equiv 0 \bmod p$. The integer $n\ne0,1$ is called a quadratic residue if the congruence has a solution $x=x_0$. Otherwise, $n$ is called a quadratic nonresidue. The Burgess upper bound of the least quadratic nonresidue claims that
\begin{equation} \label{eq9900Q.000b}
	n(p)\ll p^{\frac{1}{4\sqrt{e}}+\varepsilon},  
\end{equation} 
where $\varepsilon>0$ is a small number, see \cite {BD1963}, and \cite{BG2013} for a survey and discussion. However, in general, the least quadratic nonresidue is significantly smaller. Conditioned on the generalized RH, the upper bound satisfies
\begin{equation} \label{eq9900Q.000f}
	n(p)\ll (\log p)^2,  
\end{equation}
see \cite{AN1952}, \cite{LS2015}, et alii for fine details and explicit versions. The heuristic in {\color{red}\cite[Section 2]{MT2021}} suggests the following upper bound.
\begin{conj}\label{conj9900Q.100A} \hypertarget{conj9900Q.100A}{\normalfont (McGown-Trevino 2019)} For any $\varepsilon>0$ and a large prime $ p$, the least quadratic nonresidue has the upper bound
	$$	n(p)\ll (\log p)(\log \log p)^{1+\varepsilon}.$$
\end{conj} 

The above heuristic is slightly smaller than an earlier heuristic discussed in {\color{red}\cite[p.\;841]{EP1969}}, which leads to the conjectured upper bound
\begin{equation}
		n(p)\ll (\log p)^{1+\varepsilon}.
\end{equation}
These conjectures are phenomenal improvements of the much older Vinogradov conjecture, which claims that $n(p)\ll p^{\varepsilon}$. A closely related result due to Linnik in \cite{LY1942} proves that the number of primes $p\in [x^{\varepsilon},x]$ that fails the Vinogradov conjecture is finite, a recent proof appears in {\color{red}\cite[Corollary 5]{BR2017}}.\\

On the other extreme, there is the lower bound
\begin{equation} \label{eq9900Q.000l}
	n(p)\gg(\log p)(\log \log p),  
\end{equation}
see {\color{red}\cite[Theorem 13.5]{MH1971}}. \\

On the statistical perspective, a random integer $x<p$ is a quadratic nonresidue modulo $p$ with probability $P(x=n_p)=1/2$ and the expected value of the least quadratic nonresidue is a small constant
\begin{equation} \label{eq9900Q.000k}
	\frac{1}{\pi(x)}\sum_{p\leq x}	n_p=(1+o(1)) \sum_{n\geq 1}\frac{p_n}{2^n}=3.6746439660113287789956763\ldots,  
\end{equation} 
where $p_n$ denotes the $n$th prime in increasing order, see \cite{EP1961}. Furthermore, combining a result for Gauss quadratic sum and Weyl's theorem, it is easy to verify that both the quadratic residues and quadratic nonresidues are equidistributed on the interval $[1,p-1]$.    \\

This short note proposes a resolution of the above conjectures.

\begin{thm}\label{thm9900Q.100} \hypertarget{thm9900Q.100} Let $p$ be a large prime and let $n(p)$ denotes the least quadratic nonresidue modulo $p$. Then
	\begin{equation}\label{eq9900Q.100}
		n(p)\ll (\log p)(\log   \log p)^{1+\varepsilon},   
	\end{equation}
	where $\varepsilon>0$ is a small number.
\end{thm}

The analysis is completely different from the traditional literature on quadratic residues and quadratic nonresidues. This made possible by a new indicator function for quadratic nonresidues in finite fields introduced in \hyperlink{S9911Q-B}{Section} \ref{S9911Q-B}. The foundational and supporting materials are covered in \hyperlink{S9911Q-B}{Section} \ref{S9911Q-B} to \hyperlink{S9900Q-ET}{Section} \ref{S9900Q-ET}. 
The proof of \hyperlink{thm9900Q.100}{Theorem} \ref{thm9900Q.100} appears in  \hyperlink{S9900Q-RT}{Section} \ref{S9900Q-RT}.

\section{\textbf{Representations of the Characteristic Functions}}\label{S9911Q-B}\hypertarget{S9911Q-B}
For an odd prime $p$ the quadratic symbol modulo $p$ is defined by 
\begin{equation}\label{eq9911Q.100f}	\hypertarget{eq9911Q.100f}
	\left( \frac{n}{p}\right) 
	=\left \{
	\begin{array}{ll}
		1 & \text{ if } n \text{ is a quadratic residues},  \\[.2cm]
		-1 & \text{ if } n \text{ is a quadratic nonresidues}, \\[.2cm]
			0 & \text{ if } n \text{ is divisible by } p,  \\
	\end{array} \right .
\end{equation}
The classical characteristic functions of quadratic residues and quadratic nonresidues in the finite field $\F_p $, which are defined in terms of the quadratic symbol, have the simple formulas described below.

\begin{lem} \label{lem9911Q.200C} \hypertarget{lem9911Q.200C} If \(p\geq 2\) is a prime and \(n\in 	\mathbb{F}_p\) is a nonzero element, then
	\begin{enumerate}[font=\normalfont, label=(\roman*)]
		\item$\displaystyle 	\hypertarget{eq9911Q.100c-1}
		\varkappa _0(n)=\frac{1}{2}\left( 1+\left( \frac{n}{p}\right) \right) 
		=\left \{
		\begin{array}{ll}
			1 & \text{ if } n \text{ is n quadratic residues},  \\
			0 & \text{ if } n \text{ is n quadratic nonresidues}, \\
		\end{array} \right .$ \\[.3cm]
		
		\item$\displaystyle \varkappa (n)=\frac{1}{2}\left( 1-\left( \frac{n}{p}\right) \right) 
		=\left \{
		\begin{array}{ll}
			1 & \text{ if } n \text{ is n quadratic nonresidues},  \\
			0 & \text{ if } n \text{ is n quadratic residues}, \\
		\end{array} \right .$ 
	\end{enumerate}
	respectively.
\end{lem}	

A new representation of the characteristic function for quadratic nonresidues in the finite field $\F_p $ is introduced here.  

\begin{lem} \label{lem9911Q.200A} \hypertarget{lem9900.200A} Let \(p\geq 2\) be a prime and let \(\tau\) be a primitive root mod \(p\). If \(n\in
	\mathbb{F}_p\) is a nonzero element, then
	\begin{enumerate}[font=\normalfont, label=(\roman*)]
		\item$\displaystyle \varkappa_0 (n)=\sum _{0\leq s<p/2} \frac{1}{p}\sum _{0\leq t\leq p-1} e^ {i2\pi\frac{(\tau ^{2s}-n)t}{p}}
		=\left \{
		\begin{array}{ll}
			1 & \text{ if } n \text{ is n quadratic residues},  \\
			0 & \text{ if } n \text{ is n quadratic nonresidues}, \\
		\end{array} \right .$ \\[.3cm]
		\item$\displaystyle \varkappa (n)=\sum _{0\leq s<p/2} \frac{1}{p}\sum _{0\leq t\leq p-1} e^ {i2\pi \frac{(\tau ^{2s+1}-n)t}{p}}
		=\left \{
		\begin{array}{ll}
			1 & \text{ if } n \text{ is a quadratic nonresidues},  \\
			0 & \text{ if } n \text{ is a quadratic residues}, \\
		\end{array} \right .$ 
	\end{enumerate}
	respectively.
\end{lem}	
\begin{proof}[\textbf{Proof}] (ii) As the index \(s\geq 0\) ranges over the odd integers up to \(p-1\), the element \(\tau ^{2s+1}\in \mathbb{F}_p\) ranges over the quadratic nonresidues modulo $p$. Thus, the equation 
	\begin{equation}\label{eq9911Q.200A}
		\tau ^{2s+1}- n=0
	\end{equation} 
	has a solution if and only if the fixed element \(n\in \mathbb{F}_p\) is a quadratic nonresidue. In this case the inner sum in
\begin{equation}\label{eq9911Q.200c}
	\sum _{0\leq s<p/2} \frac{1}{p}\sum _{0\leq t\leq p-1} e^ {i2\pi\frac{(\tau ^{2s+1}-n)t}{p}}
\end{equation}	
collapses to $\sum _{0\leq t\leq p-1}1=p$. Similarly, if the fixed element \(n\in \mathbb{F}_p\) is a quadratic residue, then the equation \eqref{eq9911Q.200A} has no solution and the inner sum in \eqref{eq9911Q.200c} collapses to $\sum _{0\leq t\leq p-1}e^ {i2\pi\frac{(\tau ^{2s+1}-n)t}{p}}=0$.
\end{proof}
\section{\textbf{Finite Fourier Transform and Summation Kernels}}

\subsection{Finite Fourier Transform}
Let $f: \C \longrightarrow \C$ be a function, and let $q \in \N$ be a large integer. 

\begin{dfn}\label{dfn4400FFT.300}{\normalfont 
		The discrete Fourier transform of the function $f:\N\longrightarrow \C$ and its inverse are defined by
		\begin{equation} \label{eq4400FFT.300d}
			\hat{f}(t)=\sum_{0 \leq s\leq q-1} e^{i \pi st/q}
		\end{equation}
		and 
		\begin{equation}\label{4400FFT.300c}
			f(s)=\frac{1}{q}\sum_{0 \leq t\leq q-1}\hat{f}(m)e^{-i2\pi st/q},
		\end{equation}
		respectively.		
	}
\end{dfn} 

The finite Fourier transform and its inverse are used here to derive a summation kernel function, which is almost identical to the Dirichlet kernel, in this application $q=p$ is a prime number.

\begin{dfn} \label{dfn4400FFT.100} \hypertarget{dfn4400FFT.100}{\normalfont Let $ p$ be a prime, let $\omega=e^{i 2 \pi/p}$, and $\zeta=e^{i 2 \pi/p}$ be roots of unity. The \textit{finite summation kernel} is defined by the finite Fourier transform identity
		\begin{equation} \label{eq4400FFT.100g}
			\mathcal{K}(f(n))=\frac{1}{p} \sum_{0 \leq t\leq p-1,}  \sum_{0 \leq s\leq p-1} \omega^{t(n-s)}f(s)=f(n).\end{equation}
	} 
\end{dfn}
This simple identity is very effective in computing upper bounds of some exponential sums
\begin{equation} \label{eq4400FFT.100h}
	\sum_{ n \leq x}  f(n)= \sum_{ n \leq x}  \mathcal{K}(f(n)),
\end{equation}
where $x  < p$. 

\subsection{Summation Kernels} 
\begin{lem}   \label{lem4400SK.150A}\hypertarget{lem4400SK.150A}  Let \(p\geq 2\)  be a large primes, let $x<p-1$ and let $\omega=e^{i2 \pi/p} $ be a $p$th root of unity. If  $t \in [1, p-1]$, then,
	$$\displaystyle \left |  \sum_{n \leq x} \omega^{tn}  \right |\leq \frac{2p }{\pi t}.$$
\end{lem} 

\begin{proof}[\textbf{Proof}] Use the geometric series to compute this simple exponential sum as
	\begin{eqnarray} \label{eq4400SK.340}
		\sum_{n \leq x}\omega^{tn}
		&=& \frac{\omega^{t}-\omega^{t(p-1)}}{1-\omega^{t}} \nonumber.
	\end{eqnarray} 
	Now, observe that $\omega=e^{i2 \pi/p}$, the integers $t \in [1, p-1]$, and $d < p-1$. This data implies that $\pi t/p\ne k \pi $ with $k \in \mathbb{Z}$, so the sine function $\sin(\pi t/p)\ne 0$ is well defined. Using standard manipulations, and $z/2 \leq \sin(z) <z$ for $0<|z|<\pi/2$, the last expression becomes
	\begin{equation}
		\left |\frac{\omega^{t}-\omega^{t(x+1)}}{1-\omega^{t}} \right |\leq 	\left | \frac{2}{\sin( \pi t/ p)} \right | 
		\leq \frac{2p}{\pi t}.
	\end{equation}
\end{proof}
\begin{lem}   \label{lem4400SK.150B}\hypertarget{lem4400SK.150B}   Let \(p\geq 2\), let $x<p-1$ a and let $\omega=e^{i2 \pi/p} $ be a $p$th root of unity. Then,
	$$	\left | 	\sum_{\substack{n\leq x\\\gcd(n,p-1)=1}} \omega^{tn}  \right |\ll \frac{2p ^{1+\delta}\log p}{\pi t} , $$ 
	where $\delta>0$ is a small real number and $t \in [1, p-1]$. 
\end{lem} 

\begin{proof}[\textbf{Proof}] Set $\omega=e^{i2 \pi/p}$. Use the inclusion exclusion principle to rewrite the exponential sum as
	\begin{eqnarray} \label{eq4400SK.340d}
		\sum_{\substack{n\leq x\\\gcd(n,p-1)=1}} \omega^{tn}&=& \sum_{n \leq x} \omega^{tn}  \sum_{\substack{d \mid p-1 \\ d \mid n}}\mu(d)   \\[.3cm] 
		&=& \sum_{d \mid p-1} \mu(d) \sum_{\substack{n \leq x \\ d \mid n}} \omega^{tn}\nonumber\\[.3cm] 
		& =&\sum_{d \mid p-1} \mu(d) \sum_{m \leq (p-1)/ d} \omega^{dtm} \nonumber\\[.3cm] 
		&=& \sum_{d \mid p-1} \mu(d) \left( \frac{1-\omega^{dt(\frac{p-1}{d}+1)}}{1-\omega^{dt}}-1\right)  \nonumber\\[.3cm] 
		&=& \sum_{d \mid p-1} \mu(d)  \frac{\omega^{dt}-\omega^{dt(\frac{p-1}{d}+1)}}{1-\omega^{dt}} \nonumber,
	\end{eqnarray} 
the last 2 lines follows from a geometric summation. Taking absolute value yields
\begin{eqnarray}\label{eq4400SK.340f}
	\left|  \sum_{d \mid p-1} \mu(d) \frac{\omega^{dt}-\omega^{dtp/d}}{1-\omega^{dt}} \right| 
	&\leq&  \sum_{d \mid p-1}| \mu(d)|\cdot  \left|\frac{\omega^{dt}-\omega^{dtp/d}}{1-\omega^{dt}} \right| \\[.3cm]
	&\leq&  \sum_{d \mid p-1}  \left|\frac{\omega^{dt}-\omega^{dtp/d}}{1-\omega^{dt}} \right| \nonumber \\[.3cm]
	&\leq&  \sum_{d \mid p-1} \left | \frac{2}{\sin( \pi dt/p)} \right | \nonumber .
\end{eqnarray}
	Now, observe that the integers $t \in [1, p-1]$, and $d \leq p-1$. This data implies that $\pi dt/p\ne k \pi $ with $k \in \mathbb{Z}$. Accordingly, the sine function $\sin(\pi dt/p)\ne 0$ is well defined. In addition, For each $d\mid p-1$, the argument $dt//p/d=t/p<1$. Thus, the sine function approximation $z/2 \leq \sin(z) <z$ for $0<|z|<\pi/2$ over the subinterval $[1,p/d)$ is feasible here. Under these conditions, the last expression becomes
	\begin{eqnarray}\label{eq4400SK.340h}
\sum_{d \mid p-1} \left | \frac{2}{\sin( \pi dt/p)} \right |
&\leq& \sum_{d\mid p-1}	\left | \frac{2}{\sin( \pi dt/p/d)} \right |\\[.3cm] 
&\leq& \sum_{d\mid p-1}	\left | \frac{2}{\sin( \pi t/p)} \right |\nonumber\\[.3cm]  	 
&\leq&  \frac{2p}{\pi t}\sum_{d\mid p-1}1\nonumber\\[.3cm]
&\leq&  \frac{2p^{1+\delta}}{\pi t},\nonumber
	\end{eqnarray}
	where $\sum_{d \mid n} 1\ll n^{\delta}$, $\delta>0$. 
\end{proof}
Additional information on the order of the divisor function has a nearly explicit upper bound of the form $\sum_{d\mid n}1=n^{(\log 2+o(1)/\log\log n)}$, this appears in {\color{red}\cite[Proposition 7.12]{DL2012}}, {\color{red}\cite[Theorem 315]{HW1979}} and similar sources.

\section{\textbf{Estimates Of Exponential Sums}} \label{S9933Q}\hypertarget{S9933Q}
Exponential sums indexed by the powers of elements of nontrivial orders have applications in mathematics and cryptography. These applications have propelled the development of these exponential sums. 

\begin{thm}  {\normalfont ({\color{red}\cite[Lemma 4]{FS2002}})}   \label{thm9933Q.110A}\hypertarget{thm9933Q.110A} For integers \(a, k, N\in \mathbb{N}\), assume that \(\gcd (a,N)=c\), and that \(\gcd (k,t)=d\).
	\begin{enumerate} [font=\normalfont, label=(\roman*)]
		\item If the element \(\theta \in \mathbb{Z}_N\) is of multiplicative order \(t\geq t_0\), then 
		\begin{equation}
			\max_{1\leq a\leq p-1} \left| \sum _{ 1\leq x\leq t} e^{i2\pi a \theta ^{k x}/N} \right| <c d^{1/2}N^{1/2} .
		\end{equation} 
		\item If \(H\subset \mathbb{Z}/N \mathbb{Z}\) is a subset of cardinality \(\# H\geq N^{\delta }, \delta >0\), then 
		\begin{equation}
			\max_{\gcd (a,\varphi(N))=1}\left|  \sum _{ x\in H} e^{i2\pi a\theta ^x/N} \right| <N^{1-\delta } .
		\end{equation}
	\end{enumerate}
\end{thm}
Various upper bounds of exponential sums over subsets of elements in finite rings $\left (\mathbb{Z}/N\mathbb{Z}\right )^\times$ can be used to prove the next result. These estimates are useful in the proof of \hyperlink{lem9933RPI.220}{Lemma} \ref{lem9933RPI.220} . The reader should consult the literature, such as \cite{CC2009}, \cite{BN2004}, {\color{red}\cite[Theorem 2.1]{BJ2007}}, \cite{KS2012}, and references within the cited papers. \\

\subsection{Incomplete Exponential Sums over Consecutive Index}
The simple finite Fourier transform identity 
\begin{equation}\label{eq9933Q.110c}
	\sum_{ n \leq x}  f(n)= \sum_{ n \leq x}  \mathcal{K}(f(n)),
\end{equation}
see \hyperlink{dfn4400FFT.100}{Definition} \ref{dfn4400FFT.100}, is very effective in computing upper bounds of some exponential sums. An improved version of \hyperlink{thm9933Q.110A}{Theorem} \ref{thm9933Q.110A} and a few other applications are illustrated here.

\begin{thm}  \label{thm9933Q.110}\hypertarget{thm9933Q.110} {\normalfont (\cite{SR1973},  \cite{ML1972}) }  Let \(p\geq 2\) be a large prime, and let \(\tau \in \mathbb{F}_p\) be an element of large multiplicative order $\ord_p(\tau) \mid p-1$. Then, for any $b \in [1, p-1]$,  and $x\leq p-1$,
	\begin{equation}\label{eq9933Q.110d}
		\sum_{ n \leq x}  e^{i2\pi b \tau^{n}/p} \ll p^{1/2}  \log p.
	\end{equation}
	
\end{thm}

\begin{proof}[\textbf{Proof}]  Let $p$ be a prime, and let $f(n)=e^{i 2 \pi b\tau^{n} /p}$, where $\tau$ is a primitive root modulo $p$. Applying the finite summation kernel in \eqref{eq9933Q.110c} yields
	\begin{equation} \label{eq9933Q.110f}
		\sum_{ n \leq x}  e^{i2\pi b \tau^{n}/p}= \sum_{ n \leq x}\frac{1}{p} \sum_{0 \leq t\leq p-1,}  \sum_{1 \leq s\leq p-1} \omega^{t(n-s)}e^{i2\pi b \tau^{s}/p} .
	\end{equation}
	The term $t=0$ contributes $-x/p$, and rearranging it yield
	\begin{eqnarray} \label{eq9933Q.110h}
		\sum_{ n \leq x}  e^{i2\pi b \tau^{n}/p}
		&=&\frac{1}{p} \sum_{ n \leq x,} \sum_{1 \leq t\leq p-1,}  \sum_{1 \leq s\leq p-1} \omega^{t(n-s)}e^{i2\pi b \tau^{s}/p}-\frac{x}{p} \\[.3cm]
		&=&\frac{1}{p}  \sum_{1 \leq t\leq p-1}  \left (\sum_{1 \leq s\leq p-1} \omega^{-ts}e^{i2\pi b \tau^{s}/p} \right ) \left (\sum_{ n \leq x}\omega^{tn} \right )-\frac{x}{p}\nonumber.
	\end{eqnarray}
	Taking absolute value, and applying \hyperlink{lem440SK.150A}{Lemma} \ref{lem4400SK.150A} and \hyperlink{lem1234A.150A}{Lemma} \ref{lem1234A.150A}, yield
	\begin{eqnarray} \label{eq9933Q.110i}
		\left | \sum_{ n \leq x}  e^{i2\pi b \tau^{n}/p} \right |
		&\leq&\frac{1}{p}  \sum_{1 \leq t\leq p-1} \left | \sum_{0 \leq s\leq p-1} \omega^{-ts}e^{i2\pi b \tau^{s}/p} \right | \cdot  \left | \sum_{ n \leq x}\omega^{tn} \right |+ \frac{x}{p} \\[.3cm]
		&\ll&\frac{1}{p}  \sum_{1 \leq t\leq p-1} \left ( 2p^{1/2} \log p \right ) \cdot  \left ( \frac{2p}{\pi t} \right )+\frac{x}{p}\nonumber\\[.3cm]
		&\ll& p^{1/2} \log^2 p\nonumber .
	\end{eqnarray}
	The third line in \eqref{eq9933Q.110i} uses the estimate 
	\begin{equation} \label{eq9933Q.110j} 
		\sum_{1 \leq t\leq p-1}\frac{1}{t}\ll \log p\ll \log p.
	\end{equation} 
\end{proof}
This appears to be the best possible upper bound. The above proof generalizes the sum of resolvents method used in \cite{ML1972}. Here, it is reformulated as a finite Fourier transform method, which is applicable to a wide range of functions. A similar upper bound for composite moduli $p=m$ is also proved, [op. cit., equation (2.29)]. \\

\subsection{Incomplete Exponential Sums over Relatively Prime Index}

\begin{thm}  \label{thm9933ERP.120}\hypertarget{thm9933ERP.120} {\normalfont (\cite{SR1973}) }  Let \(p\geq 2\) be a large prime, and let \(\tau \in \mathbb{F}_p\) be an element of large multiplicative order $p-1=\ord_p(\tau)$. If $x\leq p-1$ and $s \in [1, p-1]$, then 
	\begin{equation}
		\sum_{ n \leq x}  e^{i2\pi s \tau^{n}/p} \ll p^{1/2}  \log p.
	\end{equation}
	
\end{thm}

This appears to be the best possible upper bound. A similar upper bound for composite moduli $p=m$ is also proved, [op. cit., equation (2.29)].  A simpler proof and generalization of this exponential is is provided in \cite{ML1972}.

\begin{thm}   \label{thm9933ERP.120B}\hypertarget{thm9933ERP.120B}  Let \(p\geq 2\) be a large prime, let $x\leq p$ and let $\tau $ be a primitive root modulo $p$. If $s \in [1, p-1]$, then,
	\begin{equation}
		\sum_{\substack{n\leq x\\\gcd(n,p-1)=1}} e^{i2\pi s \tau^n/p} \ll  p^{1/2+\delta} ,
	\end{equation} 
	where $\delta>0 $ is a small real number. 
\end{thm}

\begin{proof}[\textbf{Proof}] Use the inclusion exclusion principle to rewrite the exponential sum as
	\begin{eqnarray}
		\sum_{\substack{n\leq x\\\gcd(n,p-1)=1}} e^{i2\pi s \tau^n/p} &=& \sum_{n \leq p-1} e^{i2\pi s \tau^n/p}\sum_{\substack{d \mid p-1 \\ d \mid n}}\mu(d)   \\[.3cm]
		&=&\sum_{d \mid p-1} \mu(d) \sum_{\substack{n \leq p-1 \\ d \mid n}} e^{i2\pi s \tau^n/p} \nonumber \\[.3cm]
		&=&\sum_{d \mid p-1} \mu(d) \sum_{m \leq (p-1)/ d} e^{i2\pi s \tau^{dm}/p} \nonumber.
	\end{eqnarray}

	Taking absolute value, and invoking \hyperlink{thm9933ERP.120}{Theorem} \ref{thm9933ERP.120} yield 
	\begin{eqnarray}
		\left\vert  \sum_{\substack{n\leq x\\\gcd(n,p-1)=1}}e^{i2\pi s \tau^{n}/p} \right\vert 
		& \leq &\sum_{d \mid p-1} \vert \mu(d) \vert \left\vert  \sum_{m \leq (p-1)/d} e^{i2\pi s \tau^{dm}/p} \right\vert  \\[.3cm]
		&\ll&\sum_{d \mid p-1} \vert \mu(d)  \vert \left (\left (\frac{p-1}{d} \right )^{1/2} \log p \right )\nonumber \\[.3cm]
		&\ll& (p-1)^{1/2} \log(p-1) \sum_{d \mid p-1} \frac{\vert \mu(d)  \vert}{d^{1/2}} \nonumber\\[.3cm]
		&\ll&  (p-1)^{1/2} \log(p-1)\cdot p^{\delta}\nonumber\\[.3cm]
		&\ll&  p^{1/2+\delta}\nonumber .
	\end{eqnarray}
	The last inequality follows from
	\begin{equation}
		\sum_{d \mid p-1} \frac{\vert \mu(d)  \vert}{d^{1/2}} \leq   \sum_{d \mid p-1} 1\ll p^{\delta}
	\end{equation}
	for any arbitrary small number $\delta >0$, and any sufficiently large prime $p \geq 2$. 
\end{proof}

A different approach to this result appears in {\color{red}\cite[Theorem 6]{FS2000}}, and related results are given in \cite{CC2009}, \cite{FS2001}, \cite{GM2005}, and {\color{red}\cite[Theorem 1]{GK2005}}. The upper bound given in \hyperlink{thm9933Q.110}{Theorem} \ref{thm9933Q.110} seems to be optimum. A different proof, which has a weaker upper bound, appears in {\color{red}\cite[Theorem 6]{FS2000}}, and related results are given in \cite{CC2009}, \cite{FS2001}, \cite{GK2005}, and {\color{red}\cite[Theorem 1]{GK2005}}.
\subsection{Equivalent Exponential Sums over Consecutive Index} 
\label{S9933ERP}\hypertarget{S9933ERP}
For any fixed $ 0 \ne b \in \mathbb{F}_p$, the map $ \tau^n \longrightarrow b \tau^n$ is one-to-one in $\mathbb{F}_p$. Consequently, the subsets 
\begin{equation} \label{eq9933ERP.220c}
	\{ \tau^n: n\leq x \}\quad \text { and } \quad  \{ b\tau^n: n\leq x \} \subset \mathbb{F}_p
\end{equation} have the same cardinalities. As a direct consequence the exponential sums 
\begin{equation} \label{eq9933ERP.220d}
	\sum_{n\leq x} e^{i2\pi b \tau^n/p} \quad \text{ and } \quad \sum_{n\leq x} e^{i2\pi \tau^n/p},
\end{equation}
have the same upper bound up to an error term up to an error term. The result below expresses the first exponential sum in \eqref{eq9933ERP.220d} as a sum of simpler exponential sum and an error term. 

\begin{lem}   \label{lem9933ERP.220}\hypertarget{lem9933ERP.220}  Let \(p\geq 2\) be a large primes. If $\tau $ be a primitive root modulo $p$, then,
	\begin{equation} \sum_{n\leq x} e^{i2\pi b \tau^n/p} = \sum_{n\leq x} e^{i2\pi  \tau^n/p} + O(p^{1/2} \log^2 p),
	\end{equation} 
	for any $ b \in [1, p-1]$. 	
\end{lem} 
\begin{proof}[\textbf{Proof}]  For $b\ne 1$, the exponential sum has the representation 
	\begin{equation} \label{eq9933ERP.220h}
		\sum_{n\leq x} e^{\frac{i2\pi b \tau^n}{p}} 
		=\frac{1}{p} \sum_{1 \leq t\leq p-1} \left ( \sum_{1 \leq s\leq p-1} \omega^{-ts}e^{\frac{i2\pi b \tau^s}{p}}\right )\left (\sum_{n \leq x}   \omega^{tn} \right ) -\frac{\varphi(p)}{p},
	\end{equation} 
	confer equation \eqref{eq9933Q.110h} for more details. And, for $b=1$, 
	\begin{equation} \label{eq9933ERP.220i}
		\sum_{n\leq x} e^{\frac{i2\pi  \tau^n}{p}} 
		=\frac{1}{p} \sum_{1 \leq t\leq p-1} \left ( \sum_{1 \leq s\leq p-1} \omega^{-ts}e^{\frac{i2\pi  \tau^s}{p}}\right )\left (\sum_{n \leq x}   \omega^{tn} \right ) -\frac{\varphi(p)}{p}.
	\end{equation} 
	Differencing \eqref{eq9933ERP.220h} and \eqref{eq9933ERP.220i} produces 
	\begin{eqnarray} \label{eq9933ERP.220j}
		S&=&\sum_{n\leq x}e^{i2\pi b \tau^n/p} -\sum_{n\leq x} e^{i2\pi  \tau^n/p}\\ [.3cm]
		&=&     \frac{1}{p} \sum_{0 \leq t\leq p-1} \left ( \sum_{1 \leq s\leq p-1} \omega^{-ts}e^{\frac{i2\pi b \tau^s}{p}}-\sum_{1 \leq s\leq p-1} \omega^{-ts}e^{\frac{i2\pi  \tau^s}{p}}\right ) \left (\sum_{n \leq x}   \omega^{tn} \right ) \nonumber.
	\end{eqnarray}
	By \hyperlink{lem4400SK.150A}{Lemma} \ref{lem4400SK.150A}, the summation kernel is bounded by
	\begin{equation} \label{eq9933ERP.220k}
		\left |\sum_{n \leq x}   \omega^{tn}  \right | 
		\leq   \frac{2 p } {\pi t}, 
	\end{equation}
	and by \hyperlink{lem1234A.150A}{Lemma} \ref{lem1234A.150A}, the difference of two Gauss sums is bounded by
	\begin{eqnarray} \label{eq9933ERP.220l}
		G	&=&  \left | \sum_{1 \leq s\leq p-1} \omega^{-ts}e^{\frac{i2\pi b \tau^s}{p}}-\sum_{1 \leq s\leq p-1} \omega^{-ts}e^{\frac{i2\pi  \tau^s}{p}}\right  |   \\[.3cm]
		&\leq & \left |  \sum_{1 \leq s\leq p-1} \chi(s) \psi_b(s)\right |+ \left | \sum_{1 \leq s\leq p-1} \chi(s) \psi_1(s) \right | \nonumber \\[.3cm] 
		&\leq & 4 p^{1/2} \log p, \nonumber
	\end{eqnarray}
	where  $\chi(s)=e^{i \pi s t/p}$, and $ \psi_b(s)=e^{i2\pi b \tau^s/p}$. Taking absolute value in \eqref{eq9933ERP.220j} and replacing  \eqref{eq9933ERP.220k}, and  \eqref{eq9933ERP.220l}, return
	\begin{eqnarray} \label{eq9933ERP.220m}
		|S|		&=&\left|  	\sum_{1\leq n\leq p-1} e^{i2\pi b \tau^n/p} -\sum_{1\leq n\leq p-1} e^{i2\pi  \tau^n/p} \right|  \\[.3cm]
		& \leq &      \frac{1}{p} \sum_{0 \leq t\leq p-1} \left ( 4p^{1/2} \log p \right ) \cdot \left ( \frac{2 p } {t} \right )\nonumber \\[.3cm]
		&\leq & 8p^{1/2} (\log p)(\log p)( \log \log p )\nonumber\\[.3cm]
		&\leq & 8p^{1/2} \log^2 p \nonumber.
	\end{eqnarray}
\end{proof}

\subsection{Equivalent Exponential Sums over Relatively Prime Index} 
For any fixed primitive root $\tau $ and $0 \ne b \in \mathbb{F}_p$, the maps $ n \longrightarrow \tau^n$ and $ n \longrightarrow b \tau^n$ are one-to-one in $\mathbb{F}_p$. Consequently, the subsets 
\begin{equation} \label{eq9933RPI.220c}
	\{ \tau^n: n\leq x \}\quad \text { and } \quad  \{ b\tau^n: n\leq x \} \subset \mathbb{F}_p
\end{equation} have the same cardinalities. As a direct consequence the exponential sums 
\begin{equation} \label{eq9933RPI.220d}
	\sum_{\substack{n\leq x\\\gcd(n,p-1)=1}} e^{i2\pi b \tau^n/p} \quad \text{ and } \quad 	\sum_{\substack{n\leq x\\\gcd(n,p-1)=1}} e^{i2\pi \tau^n/p},
\end{equation}
have the same upper bound up to an error term up to an error term. An asymptotic formula is provided in \hyperlink{lem9933RPI.220}{Lemma} \ref{lem9933RPI.220}. The proof is based on finite Fourier transform version of the Lagrange resolvent 
\begin{equation} \label{344}
	(\omega^t,\zeta^s)=\zeta^s+\omega^{-t} \zeta^{s\tau}+\omega^{-2t} \zeta^{s\tau^{2}}+ \cdots +\omega^{-(p-1)t}\zeta^{s\tau^{p-1}}, 
\end{equation}
where $\omega=e^{i 2 \pi/p}$, $\zeta=e^{i 2 \pi/p}$, and $0 \ne s,t \in \mathbb{F}_p$. \\

The result below expresses the first exponential sum in \eqref{eq9933RPI.220d} as a sum of simpler exponential sum and an error term. 

\begin{lem}   \label{lem9933RPI.220}\hypertarget{lem9933RPI.220}  Let \(p\geq 2\) be a large primes. If $\tau $ be a primitive root modulo $p$, then,
	\begin{equation} 	\sum_{\substack{n\leq x\\\gcd(n,p-1)=1}} e^{i2\pi b \tau^n/p} = 	\sum_{\substack{n\leq x\\\gcd(n,p-1)=1}} e^{i2\pi  \tau^n/p} + O(p^{1/2+\delta}),
	\end{equation} 
	for any $ b \in [1, p-1]$ and $\delta>0$ is a small real number. 	
\end{lem} 
\begin{proof}[\textbf{Proof}]  For $b\ne 1$, the exponential sum has the representation 
	\begin{equation} \label{eq9933RPI.220h}
		\sum_{\substack{n\leq x\\\gcd(n,p-1)=1}} e^{\frac{i2\pi b \tau^n}{p}} 
		=\frac{1}{p} \sum_{1 \leq t\leq p-1} \left ( \sum_{1 \leq s\leq p-1} \omega^{-ts}e^{\frac{i2\pi b \tau^s}{p}}\right )\left (	\sum_{\substack{n\leq x\\\gcd(n,p-1)=1}}   \omega^{tn} \right ) -\frac{\varphi(p)}{p},
	\end{equation} 
	confer equation \eqref{eq9933Q.110h} for more details. And, for $b=1$, 
	\begin{equation} \label{eq9933RPI.220i}
		\sum_{\substack{n\leq x\\\gcd(n,p-1)=1}} e^{\frac{i2\pi  \tau^n}{p}} 
		=\frac{1}{p} \sum_{1 \leq t\leq p-1} \left ( \sum_{1 \leq s\leq p-1} \omega^{-ts}e^{\frac{i2\pi  \tau^s}{p}}\right )\left (	\sum_{\substack{n\leq x\\\gcd(n,p-1)=1}}   \omega^{tn} \right ) -\frac{\varphi(p)}{p}.
	\end{equation} 
	Differencing \eqref{eq9933RPI.220h} and \eqref{eq9933RPI.220i} produces 
	\begin{eqnarray} \label{eq9933RPI.220j}
		S&=&	\sum_{\substack{n\leq x\\\gcd(n,p-1)=1}}e^{i2\pi b \tau^n/p} -	\sum_{\substack{n\leq x\\\gcd(n,p-1)=1}} e^{i2\pi  \tau^n/p}\\[.3cm] 
		&=&     \frac{1}{p} \sum_{0 \leq t\leq p-1} \left ( \sum_{1 \leq s\leq p-1} \omega^{-ts}e^{\frac{i2\pi b \tau^s}{p}}-\sum_{1 \leq s\leq p-1} \omega^{-ts}e^{\frac{i2\pi  \tau^s}{p}}\right ) \left (	\sum_{\substack{n\leq x\\\gcd(n,p-1)=1}}   \omega^{tn} \right ) \nonumber.
	\end{eqnarray}
	By \hyperlink{lem4400SK.150B}{Lemma} \ref{lem4400SK.150B}, the summation kernel is bounded by
	\begin{equation} \label{eq9933RPI.220k}
		\left |	\sum_{\substack{n\leq x\\\gcd(n,p-1)=1}}   \omega^{tn}  \right | 
		\ll   \frac{2p ^{1+\delta}\log p}{\pi t}, 
	\end{equation}
	and by \hyperlink{lem1234A.150A}{Lemma} \ref{lem1234A.150A}, the difference of two Gauss sums is bounded by
	\begin{eqnarray} \label{eq9933RPI.220l}
		G	&=&  \left | \sum_{1 \leq s\leq p-1} \omega^{-ts}e^{\frac{i2\pi b \tau^s}{p}}-\sum_{1 \leq s\leq p-1} \omega^{-ts}e^{\frac{i2\pi  \tau^s}{p}}\right  |   \\[.3cm]
		&\leq & \left |  \sum_{1 \leq s\leq p-1} \chi(s) \psi_b(s)\right |+ \left | \sum_{1 \leq s\leq p-1} \chi(s) \psi_1(s) \right | \nonumber \\ [.3cm]
		&\leq & 4 p^{1/2} \log p, \nonumber
	\end{eqnarray}
	where  $\chi(s)=e^{i \pi s t/p}$, and $ \psi_b(s)=e^{i2\pi b \tau^s/p}$. Taking absolute value in \eqref{eq9933RPI.220j} and replacing \eqref{eq9933RPI.220k}, and  \eqref{eq9933RPI.220l}, return
	\begin{eqnarray} \label{eq9933RPI.220m}
		|S|		&=&\left|  	\sum_{\substack{n\leq x\\\gcd(n,p-1)=1}} e^{i2\pi b \tau^n/p} -\sum_{\substack{n\leq x\\\gcd(n,p-1)=1}} e^{i2\pi  \tau^n/p} \right|  \\[.3cm]
		& \ll &      \frac{1}{p} \sum_{0 \leq t\leq p-1} \left ( 4p^{1/2} \log p \right ) \cdot \left ( \frac{2p ^{1+\delta}\log p}{\pi t} \right ) \nonumber\\[.3cm]
		&\ll & 8p^{1/2+\delta} (\log p)(\log p)(  \log p )\nonumber\\[.3cm]
		&\ll & p^{1/2+\delta} \nonumber,
	\end{eqnarray}
	where $(\log p)^3$ is absorbed by the term $p^{\delta}$. 
\end{proof}

\subsection{FFT of Power Exponential Sum with Odd Consecutive Index} 
For any fixed $ 0 \ne b \in \mathbb{F}_p$, the map $ \tau^n \longrightarrow b \tau^n$ is one-to-one (permutation) in $\mathbb{F}_p$. Consequently, the subsets 
\begin{equation} \label{eq9933RPI.500b}
	\{ \tau^{2n+1}: n\in[1,p/2) \}\quad \text { and } \quad  \{ b\tau^{2n+1}: n\in[1,p/2) \} \subset \mathbb{F}_p
\end{equation} have the same cardinalities. As a direct consequence the exponential sums 
\begin{equation} \label{eq9933RPI.500d}
	U(b)=\sum_{1\leq n< p/2}e^{\frac{i2\pi b \tau^{2n+1}}{p}} \quad \text{ and } \quad U(1)=\sum_{1\leq n<p/2} e^{\frac{i2\pi \tau^{2n+1}}{p}},
\end{equation}
have exactly the same upper bound. An asymptotic relation for the finite Fourier transform (FFT) of the exponential sums (\ref{eq9933RPI.500d}) is provided here. 

\begin{thm}   \label{thm9933ERP.220U}\hypertarget{thm9933ERP.220U}  Let \(p\geq 2\) be a large prime. If $\tau $ be a primitive root modulo $p$ and $a=o(p)$ is not a quadratic nonresidue, then
	\begin{equation} 
		\widehat{U(a)}=	\sum_{1\leq b\leq  p-1}	 e^{ \frac{-i2\pi ab}{p}}	\sum_{1\leq n<p/2} e^{\frac{i2\pi b \tau^{2n+1}}{p}} =-  \sum_{1\leq n< p/2} e^{\frac{i2\pi  \tau^{2n+1}}{p}} + O(p^{1/2} (\log p)^2)\nonumber,
	\end{equation} 
	where the implied constant is independent of $a, b \in [1, p-1]$. 	
\end{thm} 
\begin{proof}[\textbf{Proof}] Set $x=o(p)$. For $a\in[1,x]$ and $b\in[1,p-1]$, the exponential sum has the representation 
	\begin{eqnarray} \label{eq9933RPI.500f}
		U(b)&=& \sum_{1\leq n<p/2} e^{\frac{i2\pi b \tau^{2n+1}}{p}} \\[.3cm]
		&=&\frac{1}{p} \sum_{1 \leq t\leq p-1} \left ( \sum_{1 \leq s\leq p-1} \omega^{-ts}e^{\frac{i2\pi b \tau^s}{p}}\right )\left ( \sum_{1\leq n \leq p-1}   \omega^{tn} \right ) -\frac{\varphi(p)}{p}\nonumber,
	\end{eqnarray} 
	confer equations \eqref{eq9933Q.110h} for more details. In particular, for $b=1$, 
	\begin{eqnarray} \label{eq9933RPI.500h}
		U(1)&=& 	\sum_{1\leq n<p/2}e^{\frac{i2\pi  \tau^{2n+1}}{p}} \\[.3cm]
		&=& \frac{1}{p} \sum_{1 \leq t\leq p-1} \left ( \sum_{1 \leq s\leq p-1} \omega^{-ts}e^{\frac{i2\pi a \tau^s}{p}}\right )\left ( \sum_{1\leq n <p/2}   \omega^{t(2n+1)} \right ) -\frac{\varphi(p)}{p}\nonumber,
	\end{eqnarray}
	respectively. Differencing (\ref{eq9933RPI.500f}) and (\ref{eq9933RPI.500h}) produces 
	\begin{eqnarray} \label{eq9933RPI.500i}
		U(b)-U(1)&= &	\sum_{1\leq n<p/2}e^{\frac{i2\pi b \tau^{2n+1}}{p}} -\sum_{1\leq n<p/2}e^{\frac{i2\pi  \tau^{2n+1}}{p}} \\[.3cm]
		&=&     \frac{1}{p} \sum_{1 \leq t\leq p-1} \left ( \sum_{1 \leq s\leq p-1} \omega^{-ts}e^{\frac{i2\pi  b \tau^s}{p}}-\sum_{1 \leq s\leq p-1} \omega^{-ts}e^{\frac{i2\pi  \tau^s}{p}}\right ) \left ( \sum_{1\leq n< p/2}   \omega^{t(2n+1)} \right ) \nonumber.
	\end{eqnarray}

Taking the finite Fourier transform of the difference $D(b)=V(b)-V(1)$ returns 
	
	\begin{eqnarray} \label{eq9933RPI.500j}
		\widehat{D(a)}&=&	\sum_{1\leq b\leq  p-1}	 e^{ \frac{ -i2\pi ab}{p}}\left( \sum_{1\leq n<p/2}e^{\frac{i2\pi b \tau^{2n+1}}{p}} -\sum_{1\leq n<p/2}e^{\frac{i2\pi  \tau^{2n+1}}{p}}\right)  \\[.3cm]
		&=&  \frac{1}{p} \sum_{1\leq b\leq  p-1}	 e^{ \frac{-i2\pi ab}{p}}  \sum_{1 \leq t\leq p-1} \left ( \sum_{1 \leq s\leq p-1} \omega^{-ts}e^{\frac{i2\pi  b \tau^s}{p}}-\sum_{1 \leq s\leq p-1} \omega^{-ts}e^{\frac{i2\pi  \tau^s}{p}}\right )  \nonumber \\[.3cm]
		&&\hskip3.75in \times \left ( \sum_{1\leq n< p/2}   \omega^{t(2n+1)} \right ) \nonumber\\[.3cm]
		&=&   \frac{1}{p} \sum_{1 \leq t\leq p-1} \left ( \sum_{1 \leq s\leq p-1} \omega^{-ts}\sum_{1\leq b\leq  p-1}	  e^{\frac{i2\pi b (\tau^s-a)}{p}}\right.  \nonumber \\[.3cm]
		&&\hskip 1.25in-\left .\sum_{1\leq b\leq  p-1}	 e^{-i2\pi \frac{ab}{p}}   \sum_{1 \leq s\leq p-1} \omega^{-ts}e^{\frac{i2\pi  \tau^s}{p}}\right ) \times \left ( \sum_{1\leq n< p/2}   \omega^{t(2n+1)} \right ) \nonumber.
	\end{eqnarray}
	By hypothesis, the equation $\tau^{2s+1}-a\ne0$ for any pair $a\in[1,x]$ and $s\in\{1,2,3,\ldots,p-1\}$. Thus, using the geometric sum identity $\sum_{1\leq u\leq  p-1}	 e^{i2\pi au/p}=-1$ to simplify the last expression yields
	\begin{eqnarray} \label{eq9933RPI.500l}
		\widehat{D(a)}&=&	\sum_{1\leq b\leq  p-1}	 e^{ \frac{ -i2\pi ab}{p}}\left( \sum_{1\leq n<p/2}e^{\frac{i2\pi b \tau^{2n+1}}{p}} -\sum_{1\leq n<p/2}e^{\frac{i2\pi  \tau^{2n+1}}{p}}\right)  \\[.3cm]
		&=&   \frac{1}{p} \sum_{1 \leq t\leq p-1} \left ( (-1)(-1)-(-1)  \sum_{1 \leq s\leq p-1} \omega^{-ts}e^{\frac{i2\pi  \tau^s}{p}}\right )\times \left ( \sum_{1\leq n< p/2}   \omega^{t(2n+1)} \right ) \nonumber  .
	\end{eqnarray}
	Rearranging the last equation yield
	\begin{eqnarray} \label{eq9933RPI.500k}
		\widehat{V(a)}&=&\sum_{1\leq b\leq  p-1}	 e^{-i2\pi \frac{ab}{p}}\sum_{1\leq n<p/2}e^{\frac{i2\pi b \tau^{2n+1}}{p}} \\[.3cm]
		&=& -\sum_{1\leq n<p/2}e^{\frac{i2\pi  \tau^{2n+1}}{p}}  +  \frac{1}{p} \sum_{1 \leq t\leq p-1} \left ( 1-  \sum_{1 \leq s\leq p-1} \omega^{-ts}e^{\frac{i2\pi  \tau^s}{p}}\right ) \times \left ( \sum_{1\leq n< p/2}   \omega^{t(2n+1)} \right ) \nonumber.
	\end{eqnarray}
	
	By \hyperlink{lem4400SK.150A}{Lemma} \ref{lem4400SK.150A}, the relatively prime summation kernel is bounded by
	\begin{eqnarray} \label{eq9933RPI.500m}
		\sum_{1 \leq t\leq p-1}	\Bigg |\sum_{1\leq n< p/2}   \omega^{t(2n+1)} \Bigg | 
			&\ll &  p\log p\nonumber, 
	\end{eqnarray}
and by \hyperlink{lem1234A.150A}{Lemma} \ref{lem1234A.150A}, the difference including Gauss sum is bounded by
	\begin{eqnarray} \label{eq9933RPI.500o}
		\Bigg | 1-  \sum_{1 \leq s\leq p-1} \omega^{-ts}e^{\frac{i2\pi a \tau^s}{p}}\Bigg |=	\Bigg | 1- \sum_{1 \leq s\leq p-1} \chi(s) \psi(s) \Bigg| 
		&\leq & 2 p^{1/2} \log p, 
	\end{eqnarray}
	where  $\chi(s)=e^{i \pi s t/p}$, and $ \psi(s)=e^{i2\pi a \tau^s/p}$. Taking absolute value of the remainder term
	in (\ref{eq9933RPI.500k}) and replacing (\ref{eq9933RPI.500m}), and  (\ref{eq9933RPI.500o}), return
	\begin{eqnarray} \label{eq9933RPI.500p}
	|\widehat{R(a)}|	&=&\frac{1}{p} \Bigg |\sum_{1 \leq t\leq p-1} \left ( 1-  \sum_{1 \leq s\leq p-1} \omega^{-ts}e^{\frac{i2\pi  \tau^s}{p}}\right ) \times \left ( \sum_{1\leq n< p/2}   \omega^{t(2n+1)} \right ) \Bigg | \\[.3cm]
		&=&\frac{1}{p}\sum_{1 \leq t\leq p-1}\Bigg | 1-  \sum_{1 \leq s\leq p-1} \omega^{-ts}e^{\frac{i2\pi  \tau^s}{p}} \Bigg | \cdot \Bigg |\sum_{1\leq n< p/2}   \omega^{t(2n+1)}  \Bigg | \nonumber\\[.3cm]
		&\ll &\frac{1}{p}(2 p^{1/2} \log p)\cdot \sum_{1 \leq t\leq p-1}\Bigg | \sum_{1\leq n< p/2}   \omega^{t(2n+1)}  \Bigg | \nonumber\\[.3cm]
		&\ll &\frac{1}{p}(2 p^{1/2} \log p)\cdot (p \log p) \nonumber\\[.3cm]
		&\ll & p^{1/2} (\log p)^2,\nonumber
	\end{eqnarray}
	where the implied constant is independent of $a,b\in {1,p-1}$.
\end{proof}


\begin{cor}   \label{cor9933ERP.230U}\hypertarget{cor9933ERP.230U}  Let \(p\geq 2\) be a large prime. If $\tau $ be a primitive root modulo $p$ and $0\ne a=o(p)$ is not a quadratic nonresidue, then
	\begin{equation} 
		\Bigg|\widehat{V(a)}\Bigg|=	\Bigg|\sum_{1\leq b\leq  p-1}	 e^{-i2\pi \frac{ab}{p}}	\sum_{1\leq n\leq p-1} e^{\frac{i2\pi b \tau^{2n+1}}{p}}\Bigg| = O(p^{1/2+\delta} (\log p)^2)\nonumber,
	\end{equation} 
	where $\delta>0$ is a small number and the implied constant is independent of $a$ and $ b \in [1, p-1]$. 	
\end{cor} 
\begin{proof}[\textbf{Proof}] The second line in the estimation of the upper bound in \eqref{eq9933ERP.230d} follows from \hyperlink{thm9933ERP.220U}{Theorem} \ref{thm9933ERP.220U} and the fourth line follows from \hyperlink{thm9933Q.110}{Theorem} \ref{thm9933Q.110}: 
	\begin{eqnarray} \label{eq9933ERP.230d}
		\bigg|\widehat {V(a)}\bigg |&=& 	\Bigg|\sum_{1\leq b\leq  p-1}	 e^{\frac{-i2\pi ab}{p}}	\sum_{1\leq n\leq p-1} e^{\frac{i2\pi b \tau^{2n+1}}{p}}\Bigg| \\[.3cm]
		&=&\left |- \sum_{1\leq n\leq p-1} e^{ \frac{i2\pi\tau ^n}{p}}+ O(p^{1/2} (\log p)^2 )  \right | \nonumber\\[.3cm]
		&\ll &\left |\sum_{1\leq n\leq p-1} e^{ \frac{i2\pi\tau ^n}{p}}\right |+p^{1/2} (\log p)^2  \nonumber\\[.4cm]
		&\ll&  p^{1/2} (\log p)^2\nonumber,
	\end{eqnarray}
	where the implied constant is independent of $a\ne0$. 
\end{proof}

\subsection{\textbf{Results for Gaussian Sums}}\label{lem1234A}\hypertarget{lem1234A}
Some elementary exponential sums estimates are provided in this section. 
\begin{lem}   \label{lem1234A.150A}\hypertarget{lem1234A.150A}  
	{\normalfont (Gauss sums)} Let \(p\geq 2\) be a prime, let $\chi(t)=e^{i2 \pi t/p} $ and  $\psi(t)=e^{i2\pi  \tau^t/p}$ be a pair of characters. Then, the Gaussian sum has the upper bound
	\begin{equation} \label{eq3-355}
		\left |\sum_{1 \leq t \leq p-1}    \chi(t) \psi(t) \right | \leq 2 p^{1/2} \log p.\nonumber
	\end{equation}
	
\end{lem}

\begin{lem}   \label{lem1234A.150R}\hypertarget{lem1234A.150R}   Let \(p\geq 2\) be a prime. If $\omega=e^{i 2 \pi/p}$, $\zeta=e^{i 2 \pi/p}$, and $0 \ne s,t \in \mathbb{F}_p$, then,
	the difference of two Lagrange resolvents has the upper bound
	\begin{equation} \label{eq234A.150b}
		\left | (\omega^t,\zeta^{s\tau^{dp}})-(\omega^t,\zeta^{\tau^{dp}}) \right |  \leq 2 p^{1/2} \log p. 
	\end{equation}
\end{lem} 

\begin{proof}[\textbf{Proof}] The proof for $\left | (\omega^t,\zeta^{s\tau^{dp}}) \right |  \leq  p^{1/2} \log p$ appears in \cite{ML1972}. Hence, the difference
	\begin{equation} \label{234A.150d}
		\left | (\omega^t,\zeta^{s\tau^{dp}})-(\omega^t,\zeta^{\tau^{dp}}) \right |  \leq  \left | (\omega^t,\zeta^{s\tau^{dp}}) \right | +\left |(\omega^t,\zeta^{\tau^{dp}}) \right |    \leq 2 p^{1/2} \log p.
	\end{equation}
\end{proof}

\begin{lem}   \label{lem1234A.200Q}\hypertarget{lem1234A.200Q}  
	Let \(p\geq 2\) be a prime and let $\tau\in \F_p$ be a primitive root. If $t\ne0$, then
	\begin{equation} \label{eq1234A.200Qf}
		\sum _{0\leq s<p/2}e^ {\frac{i2\pi\tau ^{2s+1}t}{p}} =\frac{w}{2}\left(\frac{(\tau t)^{-1}}{p} \right)p^{1/2},\nonumber
	\end{equation}
	where $w\ne1$ is a root of unity.
\end{lem} 
\begin{proof}[\textbf{Proof}] Rewrite the finite sum in term of the quadratic symbol $\left( a\;|\;p\right)$ in the form
	\begin{eqnarray} \label{eq1234A.200Qh}
		\sum _{0\leq s<p/2}e^ {\frac{i2\pi\tau ^{2s+1}t}{p}}&=&\frac{1}{2}\sum _{0\leq a<p}\left(1+\left(\frac{a}{p} \right)  \right) e^ {\frac{i2\pi  a\tau t}{p}}\\[.3cm]
		&=&\frac{1}{2}\sum _{0\leq a<p}\left(\frac{a}{p} \right)   e^ {\frac{i2\pi  a\tau t}{p}}\nonumber\\
		&=&\frac{1}{2}\left(\frac{(\tau t)^{-1}}{p} \right)\sum _{0\leq a<p}\left(\frac{a}{p} \right)   e^ {\frac{i2\pi  a}{p}}\nonumber\\
		&=&\frac{w}{2}\left(\frac{(\tau t)^{-1}}{p} \right)p^{1/2}\nonumber,
	\end{eqnarray}
	where $w\in\C$ is a root of unity.
\end{proof}

\section{\textbf{Fibers and Multiplicities for Quadratic Residues}} \label{S9925FMK}\hypertarget{S9925FMK}
The multiplicities of the fibers occurring in the estimate of the error term are computed in this section.
\begin{lem}  \label{lem9900Q.300S}\hypertarget{lem9900Q.300S} Let $p$ be a prime, let $ x=(\log p)^{1+\varepsilon}$ and let $\tau\in \F_p$ be a primitive root in the finite field $\F_p$.  Define the maps
	\begin{equation}\label{eq9900Q.300-m}
		\alpha(s,n)\equiv (\tau ^{2s+1}-n)\bmod p\quad \text{ and } \quad 
		\beta(u,v)\equiv uv\bmod p.
	\end{equation}	
	Then, the fibers $\alpha^{-1}(m)$ and $\beta^{-1}(m)$ of an element $0\ne m\in \F_p$
	have the cardinalities 
	\begin{equation}\label{eq9900Q.300-f}
		\#	\alpha^{-1}(m)\leq x-1\quad \text{ and }\quad \#\beta^{-1}(m)=	x
	\end{equation}
	respectively.
\end{lem}
\begin{proof}[\textbf{Proof}]Let $\mathscr{S}=\{s<p^{1-\varepsilon}\}$. Given a fixed $n\in [2,x]$, the map 
	\begin{equation}\label{eq9900Q.300-m1}
		\alpha:\mathscr{S}\times [2,x]\longrightarrow\F_p\quad  \text{ defined by }\quad  \alpha(s,n)\equiv (\tau ^{2s+1}-n)\bmod p,
	\end{equation}
	is one-to-one. This follows from the fact that the map $s\longrightarrow\tau^s \bmod p$ is a permutation the nonzero elements of the finite field $\F_p$, and the map $(s,n)\longrightarrow(\tau ^{2s+1}-n)\bmod p$ is a shifted permutation of the subset of quadratic nonresidues 
	\begin{equation}\label{eq9900Q.300-p}
		\mathscr{N}=\{\tau ^{2s+1}:s<p^{1-\varepsilon}\}\subset \F_p,
	\end{equation}
	see {\color{red}\cite[Chapter 7]{LN1997}} for more extensive details on the theory of permutation functions of finite fields. Thus, as $(s,n)\in \mathscr{S}\times [2,x] $ varies, each value $m=\alpha(s,n)$ is repeated at most $x-1$ times. Moreover, the premises no quadratic nonresidues $n\leq x=(\log p)^{1+\varepsilon}$ implies that $m=\alpha(s,n)\ne0$. This verifies that the cardinality of the fiber is
	\begin{eqnarray}\label{eq9900Q.300-f1}
		\#	\alpha^{-1}(m)&=&	\#\{(s,n)\in \mathscr{S}\times [2,x] :m\equiv (\tau ^{2s+1}-n)\bmod p\}\nonumber\\[.3cm]
		&\leq& x-1.
	\end{eqnarray}		
	Similarly, given a fixed $u\in [1,x]$, the map 
	\begin{equation}\label{eq9900Q.300-m2}
		\beta:[1,x]\times [1,p-1]\longrightarrow\F_p\quad  \text{ defined by }\quad  \beta(u,v)\equiv uv\bmod p,
	\end{equation}
	is one-to-one. Here the map $v\longrightarrow uv \bmod p$ permutes the elements of the finite field $\F_p$. Thus, as $(u,v)\in [1,x]\times [1,p-1]$ varies, each value $m=\beta(u,v)$ is repeated exactly $x$ times. This verifies that the cardinality of the fiber is
	\begin{eqnarray}\label{eq9900Q.300-f2}
		\#	\beta^{-1}(m)&=&	\#\{(u,v)\in [1,x]\times [1,p-1]:m\equiv uv\bmod p\}\nonumber\\[.2cm]
		&=&x.
	\end{eqnarray}
	
	Now each value $m=\alpha(s,n)\ne0$ (of multiplicity up to $(x-1)$ in $	\alpha^{-1}(m)$), is matched to $m=\alpha(s,n)=\beta(u,v)$ for some $(u,v)$, (of multiplicity exactly $x$ in $	\beta^{-1}(m)$). Comparing \eqref{eq9900Q.300-f1} and \eqref{eq9900Q.300-f2} prove that $\# \alpha^{-1}(m)\leq\# \beta^{-1}(m)$. 
\end{proof}
\section{\textbf{Evaluation of the Main Term}}\label{S9900Q-M}\hypertarget{S9900Q-M}
An asymptotic formula for the main term $M(x)$ is evaluated in this section.
\begin{lem} \label{lem9900M.300M}\hypertarget{lem9900M.300M} Let $\varepsilon>0$ be a small real number. If \(p\geq 2\) is a large prime and \( x=(\log p)^{1+\varepsilon}\), then
	\begin{equation}
		\sum _{2 \leq n\leq x}1 \sum _{0\leq s<p/2}\frac{1}{p}
		= \frac{x}{2}+O\left(1\right).\nonumber
	\end{equation}
\end{lem}
\begin{proof}[\textbf{Proof}] A routine calculation returns
\begin{eqnarray}
M(x)&=&\sum _{2 \leq n\leq x}1 \sum _{0\leq s<p/2}\frac{1}{p}\\[.3cm]
&=&\left(x-O(1) \right) \cdot  \frac{1}{p}\left( \frac{p}{2}+1\right) \nonumber\\[.3cm]
&=& \frac{x}{2}+O\left(1\right) \nonumber.
\end{eqnarray}
\end{proof}

\section{\textbf{Estimate For The Error Term}} \label{S9900Q-ET}\hypertarget{S9900Q-ET}
A nontrivial upper bound for the error term $E(x)$ is computed in this section. 
The error term is partitioned as $E(x)=E_{0}(x)+E_{1}(x)$. The upper bound of the first term $E_0(x)$ for $n\leq p/x$ is derived using a geometric series/sine function techniques, and the upper bound of the second term $E_1(x)$ for $p/x\leq n\leq p/2$ is derived using exponential sums techniques.
\begin{lem}  \label{lem9900P.300E}\hypertarget{lem9900P.300E} Let $\varepsilon>0$ be a small real number. Suppose \(p\geq 2\) is a large prime and \(n\leq x=(\log p)^{1+\varepsilon}\). If there is no quadratic nonresidue \(n\leq x=(\log p)^{1+\varepsilon}\), then 
	\begin{equation}\label{eq9900P.300b}
	\sum _{2 \leq n\leq x}1\sum _{0\leq s<p/2} \frac{1}{p}\sum _{1\leq t\leq  p-1} e^ {i2\pi\frac{(\tau ^{2s+1}-n)t}{p}} =O\left( (\log  p)(\log x)\right) \nonumber. 
	\end{equation} 

\end{lem}
\begin{proof}[\textbf{Proof}] The product of a point $(u,v)\in [1,x]\times [1,p/x)$ satisfies $uv<p$. This leads to the partition $[1,p/x)\cup[p/x,p/2)$ of the index $n$, which is suitable for the sine approximation $uv/p\ll\sin(\pi uv/p)\ll uv/p$ for $|uv/p|<1$ on the first subinterval $[1,p/x)$, see \eqref{eq9900P.300u1} for more details. Thus, consider the partition of the triple finite sum
	\begin{eqnarray} \label{eq9900P.300k}
		E(x)&=& \sum _{2 \leq n\leq x}
		\frac{1}{p}\sum_{s< p/2,} \sum_{ 1\leq t\leq p-1} e^{i2\pi \frac{(\tau ^{2s+1}-n)t}{p}}   \\
		&= & \sum _{2 \leq n\leq x}
		\frac{1}{p}\sum_{s< p/x,} \sum_{ 1\leq t\leq p-1} e^{i2\pi \frac{(\tau ^{2s+1}-n)t}{p}} + \sum _{2 \leq n\leq x}
		\frac{1}{p}\sum_{p/x\leq s< p/2,} \sum_{ 1\leq t\leq p-1} e^{i2\pi \frac{(\tau ^{2s+1}-n)t}{p}} \nonumber\\[.12cm]
		&=&E_{0}(x)\;+\;E_{1}(x) \nonumber.
	\end{eqnarray} 
Summing yields
	\begin{eqnarray} \label{eq9900P.300u4}
		E(x)&=& E_{0}(x)\;+\;E_{1}(x)   \\
		&\ll&  (\log x)(\log p)\;+\;  \frac{(\log p)^2}{p^{1/2}}\cdot x\nonumber\\[.12cm]
		&\ll& (\log x)(\log p)\nonumber.
	\end{eqnarray}
	This completes the estimate of the error term.
\end{proof}

\begin{lem}   \label{lem9900P.700}\hypertarget{lem9900P.700}  Let \(p\geq 2\) be a large primes and let $x=o(p)$. If $\tau $ be a primitive root modulo $p$ and there are no quadratic nonresidue $n\in[1,x]$, then,
	\begin{equation} 
E_{0}(x) = \sum _{2 \leq n\leq x}
\frac{1}{p}\sum_{s< p/x,} \sum_{ 1\leq t\leq p-1} e^{i2\pi \frac{(\tau ^{2s+1}-n)t}{p}}= O((\log x)(\log p)).
	\end{equation} 
\end{lem} 
\begin{proof}[\textbf{Proof}] To apply the geometric series/sine function techniques, the subsum $E_0(x)$ is partition as follows.
	\begin{eqnarray} \label{eq9900P.300l}
		E_0(x)&=& \sum _{2 \leq n\leq x}
		\frac{1}{p}\sum_{s< p/x,} \sum_{ 1\leq t\leq p-1} e^{i2\pi \frac{(\tau ^{2s+1}-n)t}{p}}   \\
		&= & \sum _{2 \leq n\leq x}
		\frac{1}{p}\sum_{s< p/x} \left( \sum_{ 1\leq t< p/2} e^{i2\pi \frac{(\tau ^{2s+1}-n)t}{p}}+ \sum_{ p/2\leq t\leq p-1} e^{i2\pi \frac{(\tau ^{2s+1}-n)t}{p}}\right) \nonumber\\[.12cm]
		&=&E_{0,0}(x)\;+\;E_{0,1}(x) \nonumber.
	\end{eqnarray} 
Now, a geometric series summation of the inner finite sum in the first term yields
	\begin{eqnarray} \label{eq9900P.300m}
		E_{0,0}(x)&=& \sum _{2 \leq n\leq x}
		\frac{1}{p}\sum_{s< p/x,}  \sum_{ 1\leq t< p/2} e^{i2\pi \frac{(\tau ^{2s+1}-n)t}{p}}  \\[.3cm]
		&=&   	\frac{1}{p} \sum _{2 \leq n\leq x,}\sum_{s<p/x}   \frac{e^{i2\pi (\frac{\tau ^{2s+1}-n}{p})(\frac{p}{2}+1)}-1}{1-e^{i2\pi \frac{(\tau ^{2s+1}-n)}{p}}} \nonumber\\[.3cm]
		&\leq&   	\frac{1}{p} \sum _{2 \leq n\leq x,}\sum_{s< p/x}   \frac{2}{|\sin\pi(\tau ^s-n)/p|} \nonumber,
	\end{eqnarray} 
see \cite[Chapter 23]{DH2000} for similar geometric series calculation and estimation. The last line in \eqref{eq9900P.300m} follows from the hypothesis that $u$ is not a primitive root. Specifically, $0\ne\tau^s-n\in \F_p$ for any $n \geq 1$ and any $n\leq x=(\log p)^{1+\varepsilon}$. Utilizing \hyperlink{lem9900Q.300S}{Lemma} \ref{lem9900Q.300S}, the first term has the upper bound
	\begin{eqnarray} \label{eq9900P.300u1}
		E_{0,0}(x)&=&\frac{1}{p} \sum _{2 \leq n\leq x,}\sum_{s<p/x}   \frac{2}{|\sin\pi(\tau ^{2s+1}-n)/p|}\\	[.3cm]
		&\ll&  	\frac{2}{p} \sum_{1\leq u\leq x,}\sum_{1\leq v< p/x}   \frac{1}{|\sin\pi uv/p|}\nonumber\\	[.3cm]
		&\ll&  	\frac{2}{p} \sum_{1\leq u\leq x,}\sum_{1\leq v< p}   \frac{p}{\pi uv} \nonumber\\	[.3cm]
		&\ll& 	\sum_{1\leq u\leq x}\frac{1}{u}\sum_{1\leq v< p}   \frac{1}{v} \nonumber\\	[.3cm]
		&\ll& (\log x)(\log p)\nonumber,
\end{eqnarray}
where $uv<p$ and $|\sin\pi uv/p|\ne0$ since $p\nmid uv$. Similarly, the second term has the upper bound
	\begin{eqnarray} \label{eq9900P.300v}
		E_{0,1}(x)&=& \sum_{2\leq n\leq x}
		\frac{1}{p}\sum_{s<p/x,}  \sum_{ p/2\leq t\leq p-1} e^{i2\pi \frac{(\tau ^{2s+1}-n)t}{p}}  \\[.3cm]
		&=&   	\frac{1}{p} \sum_{2\leq n\leq x,}\sum_{s<p/x}   \frac{1-e^{i2\pi (\frac{\tau ^s-n}{p})(\frac{p+1}{2})}}{1-e^{i2\pi \frac{(\tau ^s-n)}{p}}} \nonumber\\[.3cm]
		&\leq&   	\frac{1}{p} \sum_{2\leq n\leq x,}\sum_{s<p/x}  \frac{2}{|\sin\pi(\tau ^s-n)/p|} \nonumber\\[.3cm]
		&\ll&  (\log x)(\log p)\nonumber.
	\end{eqnarray}
	This is computed in the way as done in \eqref{eq9900P.300m} to \eqref{eq9900P.300u1}, mutatis mutandis. Hence, \\		
	\begin{equation}
		E_{0}(x)	=E_{0,0}(x)\;+\;E_{0,1}(x)\ll  (\log x)(\log p).
	\end{equation}
	
\end{proof}
\begin{lem}   \label{lem9900P.750}\hypertarget{lem9900P.750}  Let \(p\geq 2\) be a large primes and let $x=o(p)$. If $\tau $ be a primitive root modulo $p$ and there are no quadratic nonresidue $n\in[1,x]$, then,
	\begin{equation} 
E_{1}(x) = 	\sum _{2 \leq n\leq x}
\frac{1}{p}\sum_{p/x\leq s< p/2,} \sum_{ 1\leq t\leq p-1} e^{i2\pi \frac{(\tau ^{2s+1}-n)t}{p}}= O\left( \frac{(\log p)^2}{p^{1/2}}\cdot x\right) .
	\end{equation} 
\end{lem} 
\begin{proof}[\textbf{Proof}] The hypothesis on $n\in[1,x]$ implies that the error term is not the trivial upper bound. Next, rearrange rearrange the inner sum in the form 
\begin{eqnarray} \label{eq9900P.750i}
E_{1}(x)
&=& \sum _{2 \leq n\leq x}
\frac{1}{p}\sum_{p/x\leq s< p/2,} \sum_{ 1\leq t\leq p-1} e^{i2\pi \frac{(\tau ^{2s+1}-n)t}{p}}\\[.3cm]
&=& \frac{1}{p} \sum _{2 \leq n\leq x,} \sum _{1\leq t\leq p-1}e^ {\frac{-i2\pi nt}{p}}\sum_{p/x\leq s< p/2}  e^ {i2 \pi \frac{\tau ^{2s+1}t}{p}}  \nonumber  .
\end{eqnarray}
Taking the absolute value and an application of \hyperlink{cor9933ERP.230U}{Corollary} \ref{cor9933ERP.230U} yield
\begin{eqnarray} \label{eq9900P.750j}
E_{1}(x)&\leq&  \frac{1}{p} \sum _{2 \leq n\leq x} \Bigg| \sum _{1\leq t\leq p-1}e^ {\frac{-i2\pi nt}{p}}\sum_{p/x\leq s< p/2}  e^ {i2 \pi \frac{\tau ^{2s+1}t}{p}}\Bigg| \\[.3cm]
&\ll& \frac{(\log p)^2}{p^{1/2}}\cdot x	\nonumber.
\end{eqnarray}
This completes the estimate. 
\end{proof}
\section{\textbf{Main Result}}\label{S9900Q-RT}\hypertarget{S9900Q-RT}
Define the quadratic nonresidue counting function by
\begin{equation} \label{eq9900Q.800h}
N_p(x)	=\sum_{\substack{n\leq x\\\left( \frac{n}{p}\right)=-1}}1.
\end{equation}	 
An asymptotic formula for this function is computed below. 
\begin{proof}[{\color{blue}\normalfont\textbf{Proof of \hyperlink{thm9900Q.100}{Theorem} \ref{thm9900Q.100}}}] Let \(p>2\) be a large prime number, let $x=	(\log p)(\log   \log p)^{1+\varepsilon}$, where $\varepsilon>0$ is a small number. Suppose the least quadratic nonresidue $n>x$ and consider the sum of the characteristic function over the short interval \([2,x]\), that is, 
	\begin{equation} \label{eq9900Q.400h}
		N_p(x)=\sum _{2 \leq n\leq x}\varkappa (n)=0.
	\end{equation}
	Replacing the characteristic function, \hyperlink{lem9911Q.200A}{Lemma} \ref{lem9911Q.200A}, and expanding the nonexistence equation \eqref{eq9900Q.400h} yield
	\begin{eqnarray} \label{eq9900Q.400m}
		N_p(x)&=&\sum _{2 \leq n\leq x} \varkappa (n) \\[.3cm]
		&=&\sum _{2 \leq n\leq x}1\sum _{0\leq s<p/2} \frac{1}{p}\sum _{0\leq t\leq p-1} e^ {i2\pi\frac{(\tau ^{2s+1}-n)t}{p}} \nonumber\\[.3cm] 
		&=&\sum _{2 \leq n\leq x}1 \sum _{0\leq s<p/2} \frac{1}{p}  +\sum _{2 \leq n\leq x}1\sum _{0\leq s<p/2} \frac{1}{p}\sum _{1\leq t\leq p-1} e^ {i2\pi\frac{(\tau ^{2s+1}-n)t}{p}}\nonumber\\[.3cm] 
		&=&M(x)\; +\; E(x)\nonumber.
	\end{eqnarray} 
	
	The main term $M(x)$ is determined by $t=0$, which reduces to the exponential to \(e^{i 2\pi  st/p}=1\), it is evaluated in \hyperlink{lem9900M.300M}{Lemma} \ref{lem9900M.300M}. The error term $E(x)$ is determined by $t\ne0$, which reduces to the exponential to \(e^{i 2\pi  st/p}\ne1\), it is estimated in \hyperlink{lem9900P.300E}{Lemma} \ref{lem9900P.300E}. \\

	Substituting these estimate and replacing $x=	(\log p)(\log   \log p)^{1+\varepsilon}$ yield
	\begin{eqnarray} \label{eq9900Q.400p}
		N_p(x)&=&	\sum _{2 \leq n\leq x} \varkappa (n)		\\[.3cm]	
		&=&M(x) + E(x) \nonumber\\[.3cm]
		&=&\left[ \frac{1}{2}  x +O(1)\right] +\left[ O\left((\log p)(\log x)\right) \right] \nonumber\\[.3cm]
		&=&\frac{1}{2}  (\log p)(\log   \log p)^{1+\varepsilon}+O\left((\log p)(\log \log p)\right)  \nonumber.
	\end{eqnarray} 
	Consequently, the main term in \eqref{eq9900Q.400p} dominates the error term:
	
	\begin{eqnarray} \label{eq9900P.400v}
		N_p(x)&=&	\sum _{2 \leq n\leq x} \varkappa (n)		\\[.3cm]	
		&\gg& (\log p)(\log   \log p)^{1+\varepsilon}\left( 1+O\left(\frac{1}{(\log \log p)^{\varepsilon}}\right) \right) \nonumber\\[.3cm]
		&>&0\nonumber 
	\end{eqnarray} 
	as $p\to\infty$. Clearly, this contradicts the hypothesis \eqref{eq9900Q.400h} for all sufficiently large prime numbers $p \geq p_0$. Therefore, there exists a quadratic nonresidue \begin{equation}
		n\leq x=(\log p)(\log   \log p)^{1+\varepsilon}
	\end{equation}
	as $p\to\infty$.
\end{proof}



\end{document}